# Some Nonlinear Equations with Double Solutions: Soliton and Chaos

*Yi-Fang Chang*

*Department of Physics, Yunnan University, Kunming, 650091, China*

(E-mail: yifangchang1030@hotmail.com)

## Abstract

The fundamental characteristics of soliton and chaos in nonlinear equation are completely different. But all nonlinear equations with a soliton solution may derive chaos. While only some equations with a chaos solution have a soliton. The conditions of the two solutions are different. When some parameters are certain constants, the soliton is derived; while these parameters vary in a certain region, the bifurcation-chaos appears. It connects a chaotic control probably. The double solutions correspond possibly to the wave-particle duality in quantum theory, and connect the double solution theory of the nonlinear wave mechanics. Some nonlinear equations possess soliton and chaos, whose new meanings are discussed briefly in mathematics, physics and particle theory.

## 1. Introduction

It is well known that some nonlinear equations have the soliton solutions [1], while all nonlinear equations have the chaos solutions. The soliton and chaos possess many different characteristic: a soliton has the same shapes and velocities in a travelling process and even if through a collision, it has a definite trace which is analogous to a classical particle; the chaos is a universal phenomenon for various nonlinear systems, it describes an order and intrinsic stochastic motion which appears to be irregular and confused. Therefore, they form usually two remarkable aspects, respectively. But the both relations are being noticed increasingly. Abdullaev summarized the dynamical chaos of solitons and breathers for the sine-Gordon equation, the nonlinear Schrodinger equation and the Toda chain, etc[2]. Recently, some discussed the relations among chaos and the KdV equation [3], the perturbed sine-Gordon equation [4], the complex Ginzburg-Landau equation [5], which have the soliton solution. Warbos has even suggested an idea: chaotic solitons (chaoitons) in the conservative systems [6].

We proved that some equations have soliton solutions and chaos solutions, and their conditions are different [7]. The possible meanings of the double solutions are discussed here.

## 2. From Soliton to Chaos

2.1. The nonlinear Schrodinger equation
$$\phi_{xx} + i\phi_t + k|\phi|^2 \phi = 0, \qquad (1)$$
has a soliton solution [1]
$$\phi_s = \phi_0 \sec h[\sqrt{\frac{k}{2}}\phi_0(x - u_e t)]\exp[i(\frac{u_e}{2})(x - u_c t)], \qquad (2)$$

where the variable $\eta = x - u_e t$. Let

$$\phi = \exp[i\frac{u_e}{2}(x - u_c t)]v, \tag{3}$$

the equation (1) may become

$$\frac{dv}{d\eta} = [C + av^2 - \frac{k}{2}v^4]^{1/2}, \tag{4}$$

where $a = (u_e/2)^2 - (u_e u_c/2)$. When C=0, the soliton solution is

$$v = \sqrt{\frac{2a}{k}} \sec h^{1/2}\eta. \tag{5}$$

From this let $v = \sqrt{2a/k} \sin x$, the equation is

$$x' = \sqrt{a} \sin x, \tag{6}$$

which has the chaos solution. For a stable state whose energy is H, if $k=-b<0$, the equation will be

$$\phi'' + H\phi - b\phi^3 = 0, \tag{7}$$

whose integral is

$$\phi' = (C - H\phi^2 + \frac{b}{2}\phi^4)^{1/2}. \tag{8}$$

Let $C = H^2/2b$, so

$$\phi' = \sqrt{\frac{b}{2}}(\frac{H}{b} - \phi^2). \tag{9}$$

When $|\phi| < \sqrt{H/b}$,

$$\phi_s = \sqrt{\frac{H}{b}} th(\sqrt{\frac{b}{2}}\eta + C_0). \tag{10}$$

It is the simplest soliton with a bell shape. Using a substitution $\phi = Hx/\sqrt{2b}$ for Eq. (9), and it become a difference equation

$$X_{n+1} = 1 - \frac{H}{2}X_n^2. \tag{11}$$

It is a known equation, which has the chaos solution, and its parameter determined the bifurcation-chaos is $\mu = H/2$. Moreover, this equation may include the Higgs equation and the Ginzburg-Landau equation.

2.2. The Dirac equation has shown the existence of a nondegenerate, isolated, zero-energy, c-number solution. Its solutions may be monopoles, dyons and solitons [8,9,10]. The nonlinear Dirac equation is

$$\gamma_\mu \partial_\mu \psi + m\psi - l_0^2 \psi(\psi^+ \psi) = 0. \tag{12}$$

It is the Heisenberg unified equation [11] when $m=0$. The probability density $\rho = \psi^+ \psi = \overline{\psi}\gamma_4 \psi$,

$$\frac{\partial \rho}{\partial x_\mu} = \frac{\partial \psi^+}{\partial x_\mu}\psi + \psi^+ \frac{\partial \psi}{\partial x_\mu} = [-\gamma_\mu l_0^2 \psi^+(\psi^+\psi) + \gamma_\mu m\psi^+]\psi + \psi^+[\gamma_\mu l_0^2 \psi(\psi^+\psi) - \gamma_\mu m\psi]$$

and $\psi^+\psi = 1 - \psi\psi^+$, so

$$\frac{\partial \rho}{\partial x_\mu} = 2\gamma_\mu l_0^2(\psi^+\psi)(\psi^+\psi) - \gamma_\mu l_0^2(\psi^+\psi) = \gamma_\mu l_0^2 \rho(2\rho - 1). \tag{13}$$

Let $\eta = \alpha(\gamma_\alpha x_\alpha - u\gamma_0 t)$, the equation is

$$\frac{d\rho}{d\eta} = l_0^2 \rho(2\rho - 1), \tag{14}$$

whose solution is

$$\rho = \frac{1}{2 - \exp(-l_0^2 \eta + c)}. \tag{15}$$

It is analogous to a soliton since $\rho = (2 - e^c)^{-1} (\eta = 0)$ and $\rho = 1/2 (\eta \to \infty)$. Using a substitution $\rho = (2 - l_0^2 x)/8$ for Eq. (14), then the corresponding difference equation is

$$X_{n+1} = 1 - \frac{1}{4} l_0^4 X_n^2, \tag{16}$$

which has the chaos solution and the parameter $\mu = l_0^4 / 4$.

2.3. For the Korteweg-de Vries equation

$$\psi_t + \sigma \psi \psi_x + \psi_{xxx} = 0, \tag{17}$$

let $\eta = x - ut$, using two order integrals, then

$$\psi' = (-\frac{1}{3}\sigma\psi^3 + u\psi^2 + C_0\psi + C_1)^{1/2}. \tag{18}$$

For the soliton solution, the integral constants should be $C_0 = C_1 = 0$, so Eq. (18) is

$$\psi' = \psi(u - \frac{1}{3}\sigma\psi)^{1/2}, \tag{19}$$

whose soliton solution is

$$\psi = \frac{3u}{\sigma} \sec h^2 (\frac{\sqrt{u}}{2}\eta). \tag{20}$$

Using a substitution $\psi = 3u[1 - (u/4)(-x)^2]/\sigma$, the difference equation is

$$X_{n+1} = 1 - \frac{1}{4} u X_n^2. \tag{21}$$

In a $0 \leq u \leq 8$ region, the values of bifurcation-chaos are $u=3,5,..., 5.6046207...$

2.4. For the cubic Klein-Gordon equation

$$\square \varphi - m^2 \varphi + a\varphi^3 = 0, \tag{22}$$

let $\eta = (x - ut)/\sqrt{1 - u^2}$, so

$$\frac{d\varphi}{d\eta} = (-\frac{1}{2} a\varphi^4 + m^2 \varphi^2 + C)^{1/2}, \tag{23}$$

If $C=0, a>0$,

$$\varphi = \sqrt{\frac{2}{a}} m \sec h(m\eta + C_0). \tag{24}$$

It is the simplest soliton with a kink shape. Moreover,

$$\frac{d\varphi}{d\eta} = m\varphi(1 - \frac{a}{2m^2}\varphi^2)^{1/2} \tag{25}$$

is the same with Eq. (4), so it has the chaos solution. Further, all nonlinear equations have chaos.

**3. From Chaos to Soliton**

3.1. The simplest difference equation with the chaos solution is

$$X_{n+1} = 1 - \mu X_n^2. \tag{26}$$

It may correspond to a differential equation of first order

$$x' = 1 - \mu x^2, \tag{27}$$

and a partial differential equation of second order

$$\phi_{xx} - \phi_{tt} + a\phi - b\phi^3 = 0. \tag{28}$$

It becomes to an ordinary differential equation by a way on soliton solution, i.e., Eq. (7). When $|x| < 1/\sqrt{\mu}$ for Eq. (26),

$$x = \frac{1}{\sqrt{\mu}} th(\mu\eta + C) \tag{29}$$

is namely a soliton solution. A bifurcation-chaos region $2 \geq \mu \geq 0, x \in [-1,1]$ corresponds to $\infty \geq \mu^{-1/2} \geq 1/\sqrt{2} = 0.7071...$ For single stable solution $0.75 \geq \mu_0 \geq 0$, i.e., $\infty \geq \mu_0^{-1/2} \geq 1.154$, so the condition on $|x| < \mu^{-1/2}$ is satisfied necessarily in the region, the soliton can exist. While for two-branch region, $1.25 \geq \mu_1 \geq 0.75$, i.e., $1.154 \geq \mu_1^{-1/2} \geq 0.894$; for a region from four-braches to chaos, $1.401152 \geq \mu \geq 1.25$, i.e., $0.894 \geq \mu^{-1/2} \geq 0.845$. Since $|x| \leq 1$, the necessary condition in which the soliton appears is $\mu \leq 1$, it corresponds to the region of single solution and a part of two-branch region. For the rest $|x| < \mu^{-1/2}$ does not hold generally.

3.2.The logistic equation

$$\frac{dF}{dt} = F(\alpha E - F) \tag{30}$$

corresponds to a difference equation

$$X_{n+1} = 1 - \frac{1}{4}\alpha^2 E^2 X_n^2, \tag{31}$$

whose parameter is $\mu = (\alpha E)^2 / 4$. In the $0 \leq \mu \leq 2$ region, two branches appear for $\alpha E = \sqrt{3}$, four branches appear for $\alpha E = \sqrt{5}$, etc., there is the chaos for $\alpha E = (5.6046207)^{1/2} \approx 2.37$. The equation (30) has the solution

$$F = \frac{\alpha E}{1 + C\exp(-\alpha Et)}. \tag{32}$$

When $t \geq 0$, Eq. (32) is analogous to a soliton since $F = \alpha E / (1+C)$ for $t=0$ and $F = \alpha E$ for $t \to \infty$. It shows that the state will reach to stable at last as time increases continuously.

3.3.The difference equation with a chaos solution

$$X_{n+1} = \lambda \sin(\pi X_n) \tag{33}$$

corresponds to a partial differential equation of second order

$$\phi_{xx} - \phi_{tt} = \frac{1}{2}\pi\lambda^2 \sin(2\pi\phi), \tag{34}$$

namely, the sine-Gordon equation. It has the soliton solution

$$\phi = \frac{2}{\pi} tg^{-1}[\exp(\pm\pi\lambda\frac{x-ut}{\sqrt{1-u^2}})]. \tag{35}$$

Only some chaos equations have the soliton solutions.

**4.Discussion**

Further, we discuss some possible meanings of the double solutions possessed by these equations briefly.

In the mathematical aspect, when some parameters are a certain constant, the soliton is derived; while these parameters vary in a certain region, the bifurcation-chaos appears. Therefore, the former corresponds to a stable state, and the latter is a changeable process.

In the physical aspect, Szebehely and McKenzie discussed that the three-body problem in gravitational field possesses chaotic behaviors [12]. We proved that the gravitational wave is a type of nonlinear wave, and should be different to electromagnetic wave and have new characteristics, for example, as solitons [13]. Perhaps, the double solutions are two different states. These parameters are the order parameters. These states often depend on the integral constants, the boundary conditions and the initial conditions. It explains again that the solutions of the nonlinear equations depend on the initial values sensitively. It connects the chaos control by a method of parameter-control. When we control the order parameter in the nonlinear system, chaos appear, disappear, synchronize [14], even a determinational soliton is produced, for different parameteric values. For example, the soliton can be derived in a propagation of shallow water waves, but if the flow rate reaches a certain value, there will form the turbulence. Moreover, the soliton solution corresponds to particle even it may be a degenerate doublet with Fermi number $\pm(1/2)$ [7,8]. The chaos solution seems to correspond to the field, including the stochastic field. It will probably connect the double solution theory of the de Broglie-Bohm nonlinear wave mechanics. In this case the wave-particle duality is a wave-particle synthesis, where the particle is described by the mobile singularity of soliton of the wave equation [15].

The double solutions show a simultaneous existence on determinism and probabilism quantitatively from an aspect in some nonlinear systems.

The soliton equations and the chaos equations have the widely applied domains, in which above double solutions will show many meaning results or a good deal of enlightenment.